\documentclass[11pt,a4paper]{amsart}

\usepackage{geometry}
\pdfoutput=1
\usepackage{amsmath,amsfonts,amssymb,amsthm}
\usepackage{xcolor}
\usepackage{todonotes}
\usepackage{hyperref}
\hypersetup{
	pagebackref  = true,
	colorlinks   = true,
	pdfauthor    = {Pavlo Yatsyna},
  pdftitle     = {A lower bound for the rank of a universal quadratic form with integer coefficients in a totally real number field},
	pdfkeywords  = {},
	pdfcreator   = {Pdflatex},
	pdfpagemode  = UseNone,
	pdfstartview = FitH
}
\usepackage{pdfpages}
\usepackage[normalem]{ulem}
\usepackage{url}
\usepackage{enumerate}
\usepackage{verbatim}
\usepackage{mathrsfs}
\usepackage{csquotes}
\usepackage{mathtools}
\usepackage{lipsum}
\numberwithin{equation}{section}

\theoremstyle{plain}
\newtheorem{theorem}{Theorem}
\newtheorem{lemma}[theorem]{Lemma}

\newtheorem{corollary}[theorem]{Corollary}
\newtheorem{proposition}[theorem]{Proposition}
\numberwithin{equation}{section}

\theoremstyle{definition}
\newtheorem{definition}[theorem]{Definition}

\newtheorem{example}{Example}

\usepackage{tikz}

\DeclarePairedDelimiterX\set[2]\lbrace\rbrace{\,#1\mathclose{}:\mathopen{}#2\,}

\title[]{A lower bound for the rank of a universal quadratic form with integer coefficients in a totally real number field}

\author{Pavlo Yatsyna}
\email{pvyatsyna@gmail.com}

\begin{document}

\begin{abstract}
	We show that if $K$ is a monogenic, primitive, totally real number field, that contains units of every signature, then there exists a lower bound for the rank of integer universal quadratic forms defined over $K$. In particular, we extend the work of Blomer and Kala, to show that there exist infinitely many totally real cubic number fields that do not have a universal quadratic form of a given rank defined over them. For the real quadratic number fields with a unit of negative norm, we show that the minimal rank of a universal quadratic form goes to infinity as the discriminant of the number field grows. These results follow from the study of interlacing polynomials. Specifically, we show that there are only finitely many irreducible monic polynomials related to primitive number fields of a given degree, that have a bounded number of interlacing polynomials.
\end{abstract}
\maketitle
\section{Introduction}\label{sec:1}

Theorem 290 (\cite{5}) tells us when a given positive definite quadratic form over rational integers represents all natural numbers. Relatively little is known about universal quadratic forms with coefficients in the ring of integers $R$ of a totally real number field $K$. G\"otzky (\cite{14}), Maass (\cite{22}) and Siegel (\cite{28}) showed that if all of the totally positive integers in $K$ can be represented as the sum of squares, then $K=\mathbb{Q}$ or $\mathbb{Q}(\sqrt{5})$. Chan, Kim, and Raghavan (\cite{7}) proved that one finds universal ternary quadratic forms only over $\mathbb{Q}(\sqrt{2})$, $\mathbb{Q}(\sqrt{3})$ and $\mathbb{Q}(\sqrt{5})$. On another hand, B. M. Kim (\cite{20}) gave an explicit construction of an octonary diagonal universal quadratic form for infinitely many real quadratic number fields. Blomer and Kala in \cite{6} (and Kala in \cite{18}) proved that for a given $m$ there are infinitely many real quadratic number fields that do not admit universal quadratic forms with $m$ coefficients. This result was extended to multiquadratic fields by Kala and Svoboda in \cite{19}. For a totally real number field $K$ of an odd degree, it was shown in \cite{10}, by Earnest and Khosravani, that there are at most finitely many inequivalent quaternary universal quadratic forms over $K$. Kitaoka suggested (\cite{21}) that there may only exist finitely many totally real number fields over which there exists a universal ternary quadratic form. 

In this paper, we investigate universal quadratic forms defined over the ring of integers $R$ of a totally real number field $K$ by studying interlacing polynomials. Specifically, we establish a correspondence between elements of the codifferent of $R$ over $K$ and interlacing polynomials. For primitive extensions we prove the following result:

\begin{theorem}\label{Thm:1}
	Let $d,r \in \mathbb{N}$. Up to $\mathbb{Z}$-equivalence, there are only finitely many irreducible monic polynomials $f\in \mathbb{Z}[x]$ of degree $d$ such that $\mathbb{Q}[x]/(f)$ is a primitive field and $f$ is interlaced by at most $r$ monic integer polynomials.
\end{theorem}

We show that there exists an infinite family of non-interlacing polynomials:

\begin{theorem}\label{Thm:2}
	Let $d\in\mathbb{N}$ such that $d$ is square-free, $d$ is not a prime number or twice a prime number, $d>20$ and $d\not=30$. Then the minimal polynomial of $\zeta_{d}^{}+\zeta_{d}^{-1}$ is non-interlacing. 
\end{theorem}

But, such polynomials do not exist for the first few degrees:

\begin{theorem}\label{Thm:3}
	There does not exist a non-interlacing irreducible integer polynomial of degree $2,3,4,5$ or $7$.
\end{theorem}

Let us say that a totally real number field satisfies Condition (A) if it is monogenic and has units of every signature. In the spirit of the work of Blomer and Kala, for primitive number fields we show:

\begin{theorem}\label{Thm:4}
	Let $d,r \in \mathbb{N}$. There are only finitely many totally real primitive number fields of degree $d$ that satisfy Condition (A) and have an integer universal quadratic form of rank $r$ defined over them. 
\end{theorem}

As an application of the theorem above, we have:

\begin{theorem}\label{Thm:5}
	For any given $r\in \mathbb{Z}$, there exist infinitely many totally real quadratic and cubic number fields that do not have a universal quadratic form of rank $r$ defined over them.
\end{theorem}

We also show that the minimal rank of a universal quadratic form over quadratic fields grows at the order of the fourth root of the discriminant (Corollary \ref{c:1-5}). This was previously shown in \cite[Prop. 5]{6} conditionally on the Riemann hypothesis.

The basic idea of this paper stems from the observation that if the codifferent of $R$ is a principal ideal generated by a totally positive number, i.e. $R^{\vee}=\gamma R$, then a universal quadratic form $Q$ over $R$ can be scaled to represent all the totally positive elements in the codifferent. Furthermore, $\mathrm{tr}_{K/\mathbb{Q}}(\gamma Q)$ becomes a quadratic form over $\mathbb{Z}$, where all the minimal vectors of $\mathrm{tr}_{K/\mathbb{Q}}(\gamma Q)$ relate to the totally positive numbers of the minimal trace in the codifferent represented by $\gamma Q$. On the other hand, totally positive elements of trace one correspond to lattice points in certain convex sets, and in particular, to interlacing polynomials. By showing that the number of interlacing polynomials grows to infinity for a given degree, we are able to show that the rank of universal quadratic forms over $R$ has to grow also, else we would have quadratic forms with the number of minimal vectors surpassing the kissing number, which is impossible. 

In the first section, after the preliminaries, we focus on the interlacing polynomials and show the correspondence between lattice points in certain convex sets and totally positive elements of trace one in the codifferent of $R$. Theorem \ref{Thm:1} appears in this section. In the following section, we tackle Theorems \ref{Thm:2} and \ref{Thm:3} and look at some of the examples of non-interlacing polynomials. Finally, the last section consists of results relating to universal quadratic forms. Appropriately, it contains proofs of Theorem \ref{Thm:4} and Theorem \ref{Thm:5}.

\section{Acknowledgements}
I would like to thank James McKee for his valuable comments on the first manuscript. This research was supported through the programme "Oberwolfach Leibniz Fellows" by the Mathematisches Forschungsinstitut Oberwolfach in 2017. 

\section{Preliminaries}

Throughout the paper, $K$ is a totally real number field with the ring of integers $R$. A number field is \textit{monogenic} if its ring of integers has integer power basis, i.e. $R=\mathbb{Z}[\alpha]$. We say that $K$ is \textit{primitive} if there is no proper subfield of $K$ other than $\mathbb{Q}$. We let $\sigma_1,\ldots,\sigma_d$ to be the distinct embedding of $K$ into $\mathbb{R}$. For $\alpha\in K$, we write $\mathrm{N}_{K/\mathbb{Q}}(\alpha)=\prod^d_{i=1}\sigma_i(\alpha)$ and $\mathrm{tr}_{K/\mathbb{Q}}(\alpha)=\sum^d_{i=1}\sigma_i(\alpha)$. For an ideal $I$ in $R$, $\mathrm{N}(I)$ is the absolute norm of $I$. We assume that all the polynomials are integer, monic, separable, and have only real roots unless stated otherwise. Let $f\in \mathbb{Z}[x]$ be an irreducible polynomial such that $K\cong \mathbb{Q}[x]/(f)$. Let $\alpha\in \mathbb{R}$ be a root of $f$, then $\mathbb{Z}[\alpha]$ is an \textit{order} in $K$. We denote the \textit{dual} of it by $\mathbb{Z}[\alpha]^{\vee}=\{\beta\in K\mid \mathrm{tr}_{K/\mathbb{Q}}(\beta\mathbb{Z}[\alpha])\subset\mathbb{Z}\}.$ For the ring of integers $R$ of $K$, its dual is the \textit{codifferent} ideal, $R^{\vee}$.
An element $\alpha\in K$ is \textit{totally positive}, and denoted by $\alpha\gg 0$, if $\sigma_i(\alpha)>0$ for all the embedding of $K$ in $\mathbb{R}$. For the set of all the totally positive numbers and integers, we write $K_{+}$ and $R_{+}$, respectively. 
We denote by $\mathrm{Disc}(f)$ the discriminant of a polynomial $f$, and by $\mathrm{Disc}(K)$ the discriminant of a number field $K$. Let $g\in \mathbb{R}[x]$, we write $\mathrm{Res}(f,g)$ for the resultant of $f$ and $g$.  

\begin{definition}\label{def:1-1}
	Let $f=\prod^d_{i=1}(x-\alpha_i),g=\prod^{d-1}_{i=1}(x-\beta_i)\in\mathbb{Z}[x]$, where $d\ge 2$. Then $g$ \textit{interlaces} $f$, or $f$ is \textit{interlaced} (by $g$), if 
	\[
	\alpha_1<\beta_1<\alpha_2<\ldots<\alpha_{d-1}<\beta_{d-1}<\alpha_d.
	\]
\end{definition}
Let us denote by 
\[\mathrm{Span}(\alpha)=\max_{i,j}|\alpha_i-\alpha_j|\]
the \textit{span} of a totally real algebraic number $\alpha$, where $\alpha=\alpha_1,\ldots,\alpha_d$ are all of the conjugates of $\alpha$. Note that $\mathrm{Span}(m)=0$ if and only if $m \in \mathbb{Q}$. And $\mathrm{Span}(\alpha)=\mathrm{Span}(\alpha+k)$ for all $k \in \mathbb{Z}$. 
We say that $f,g \in \mathbb{Z}[x]$ are \textit{$\mathbb{Z}$-equivalent} if $f(x)=g(x+n)$ for $n\in \mathbb{Z}$.

We say that a quadratic form $Q(v_1,\ldots,v_r)=\sum_{1\le i\le j \le r}a_{ij}v_iv_j$ is \textit{integral} (or over $R$) if $a_{ij}\in R$ for all $i,j$. Furthermore, those integral quadratic forms that have $a_{ij}\in2R$ whenever $i\not=j$ are called \textit{classical}. If $Q(v)\gg0$ for all $v\in R^r$, $v\not=0$, then $Q$ is said to be \textit{totally positive definite}. An $R$-lattice $(L,B_Q)$ is a pair, where $L$ is a free $R$-module of finite rank, with a non-degenerate symmetric bilinear form $B_Q:L\times L\rightarrow R$ and the associated quadratic form $Q(v)=B_Q(v,v)$ for $v \in L$. Thus, an $R$-lattice naturally corresponds to a classical quadratic form over $R$. An $R$-lattice $(L,B_Q)$ is positive definite if and only if $Q$ is positive definite. Of special importance to us will be \textit{ideal lattices}, i.e. those lattices that can be represented in the form $(I,\mathrm{tr}_{K/\mathbb{Q}}(\gamma xy))$, where $I$ is a fractional ideal in $K$ and $\gamma \in K$ (see \cite{3,4}).

We say that a quadratic form $Q$ over $R$ \textit{represents} $\alpha$ if there exits $v\in R^r$ such that $Q(v)=\alpha$. A positive definite quadratic form is \textit{universal} over $R$ if it represents all elements in $R_+$. As a convention, we always assume that the universal quadratic form is positive definite. We say that a lattice is universal if the corresponding quadratic form is universal. For a positive definite $\mathbb{Z}$-lattice $L$, we write 
\[\min(L,B_Q)=\min\{Q(v)\mid v\in L\}\]
for the minimum of $L$, and 
\[\mathcal{M}(L,B_Q)=\{v\in L\mid Q(v)=\min(L,B_Q)\}\]
for the set of the minimal vectors of $L$ (we may write $\mathcal{M}(L)$ if the bilinear form is clear from the context). Let us write 
\[
\tau_d=\max_{(L,B_Q)}|\mathcal{M}(L,B_Q)|
\]
for the \textit{kissing number} for lattices in dimension $d$, where $L$ ranges through all the positive definite $\mathbb{Z}$-lattices of rank $d$. We shall write 
\[
\mathcal{M}(K)=\{\alpha \in R^{\vee}_+\mid \mathrm{tr}_{K/\mathbb{Q}}(\alpha)=\min_{\gamma \in R^{\vee}_+}\mathrm{tr}_{K/\mathbb{Q}}(\gamma )\}.
\]
For lattices $L_1$ and $L_2$ we write $L=L_1\perp L_2$, if $L=L_1\oplus L_2$ and $B_Q(l_1,l_2)=0$ for all $l_1\in L_1$, $l_2\in L_2$. For $\alpha\in K$, we define a symmetric bilinear form scaled by $\alpha$ 
\[
B_Q^{\alpha}(v,v)=\alpha B_Q(v,v)
\]
and a quadratic form scaled by $\alpha$
\[
Q^\alpha(v)=\alpha Q(v).
\]	
The following identity clearly holds $B_Q^{\alpha}=B_{Q^{\alpha}}$.

Let $\mathrm{Mat}(d,K)$ denote the set of all square $d\times d$ matrices, and $\mathrm{diag}(a_1,\ldots,a_d)$ is a $d\times d$ diagonal matrix with elements $a_i$ on the diagonal. For a matrix $A\in \mathrm{Mat}(d,K)$, $A^t$ and $\mathrm{det}(A)$ denotes its transpose and determinant, respectively. For a given set $S$, we write $|S|$ for the cardinality of $S$. 

\begin{definition}
	Let $(L,B_Q)$ be an $\mathbb{Z}$-lattice in $\mathbb{R}^d$, and let $C\subset\mathbb{R}^d$ be a convex set. Then the \textit{width} of $C$ is \[
	w(C,L)=\min\{\max_{w\in C}B_Q(v,w)-\min_{ w\in C}B_Q(v,w) \mid v \in L\setminus\{0\}\},
	\]
	where $B_Q$ is defined on $L\times L$, extended to $L\times\mathbb{R}^d$.
\end{definition}

Finally, $\mathrm{conv}\{b_1,\ldots,b_k\}$ denotes the \textit{convex hull} of elements $b_1,\ldots,b_k\in \mathbb{R}^d$.

\section{Interlacing polynomials}

Our definition of interlacing (see Definition \ref{def:1-1}) requires that both polynomials are integer, and throughout the paper, the term interlaced will refer to integer polynomials unless stated otherwise. However, interlacing is defined analogously for real polynomials, and for interlacing with real polynomials we have:

\begin{proposition}\textup{\cite{16}}\label{prop:1-1}
	The set of all the polynomials $g\in\mathbb{R}[x]$ that interlace $f=\prod^d_{i=1}(x-\alpha_i)$ forms a convex set with vertices $\prod^{d}_{j\not=i}(x-\alpha_j)$ for $i=1,\ldots,d.$
\end{proposition}

We can associate any monic polynomial $g$ of degree $d-1$ with a vector in $\mathbb{R}^d$ as follows:
\begin{align}
\{g\in\mathbb{R}[x]\mid \deg(g)=d-1\}&\longrightarrow \mathbb{R}^d\label{eq:1}\\
x^{d-1}+\sum^{d-1}_{i=1}a_ix^{d-1-i}&\mapsto (1,a_1,\ldots,a_{d-1})^t.\nonumber
\end{align} 
Thus finding interlacing polynomials of $f$ is equivalent to finding integer points in the convex set \[
\mathcal{K}(f)=\mathrm{conv}\{(1,b^{(i)}_1,\ldots,b^{(i)}_{d-1})^t|i=1,\ldots,d\},
\]
where $\prod^{d}_{j\not=i}(x-\alpha_j)=x^{d-1}+\sum^{d-1}_{j=1}b^{(i)}_jx^{d-1-j}$.

\begin{lemma}\label{lem:1-1}
	Let $A=A_{\mathcal{K}(f)}\in\mathrm{Mat}(d,\mathbb{R})$, such that its columns are the vertices of $\mathcal{K}(f)$. Then $A^{-1}_{ij}=\dfrac{\alpha^{n-j}_i}{f'(\alpha_i)}$.
	\begin{proof}
		By a direct computation it follows that 
		\begin{align*}
		\delta_{ij}&=\dfrac{1}{f'(\alpha_i)}\prod^d_{\substack{k=1\\k\not= j}}(\alpha_i-\alpha_k)\\
		&=\dfrac{1}{f'(\alpha_i)}\sum^d_{k=1}\alpha^{d-k}_i \sum_{\substack{1\leq m_1<\ldots<m_{k-1}\leq d\\m_i\not=j}}(-1)^{k-1}\alpha_{m_1}\ldots\alpha_{m_{k-1}}\\
		&=\sum^d_{k=1}\dfrac{\alpha^{d-k}_i }{f'(\alpha_i)}A_{kj}\\
		&=(A^{-1}A)_{ij}.\qedhere
		\end{align*}
	\end{proof}
\end{lemma}

For $\mathcal{K}(f)\subset \mathbb{R}^{d}$ we define the set \[
\Lambda_{\mathcal{K}(f)}=\{\lambda \in \mathbb{R}^d_+|A_{\mathcal{K}(f)}\lambda\in \mathbb{Z}^{d}, \;\sum^d_{i=1}\lambda_i=1\}.
\]
Clearly, $\mathcal{K}(f)\cap\mathbb{Z}^{d}\not=\emptyset$ if and only if $\Lambda_{\mathcal{K}(f)}\neq \emptyset$. And for any polynomial, $\mathcal{K}(f)$ being compact implies that $|\mathcal{K}(f)\cap \mathbb{Z}^d|$ is a finite set. 	

\begin{proposition}\label{prop:1-2}
	Let $f\in \mathbb{Z}[x]$ and ${\lambda}\in \Lambda_{\mathcal{K}(f)}$, and let $g\in\mathbb{Z}[x]$ be the polynomial that interlaces $f$ corresponding to $\lambda$. Then $\left|\prod^d_{i=1}\lambda_i\right|=\left|\dfrac{\mathrm{Res}(f,g)}{\mathrm{Disc}(f)}\right|$.
	\begin{proof}
		Let $A=A_{\mathcal{K}(f)}$ be the matrix as defined in the lemma above, and let ${\lambda}\in \Lambda_{\mathcal{K}(f)}$. Now, 
		\[
		A\lambda=b=(1,b_1,\ldots,b_{d-1})^t \in \mathbb{Z}^{d},
		\] 
		and 
		\begin{align*}
		g&=x^{d-1}+\sum^{d-2}_{i=1}b_ix^{d-1-i}\\
		&=\sum^d_{i=1}\lambda_i\prod_{j\not =i}(x-\alpha_j)
		\end{align*}
		interlaces $f$.  As $\lambda=A^{-1}b$, by applying Lemma \ref{lem:1-1} it follows that
		\begin{align}
		\lambda_i&=\sum^d_{j=1}A^{-1}_{ij}b_j\nonumber\\
		&=\dfrac{1}{f'(\alpha_i)}\sum^d_{j=1}\alpha^{d-j}_ib_j\nonumber\\
		&=\dfrac{g(\alpha_i)}{f'(\alpha_i)}.\label{eq:2}
		\end{align}
		Therefore, 
		\[\prod^d_{i=1}\lambda_i=\prod^d_{i=1}\dfrac{g(\alpha_i)}{f'(\alpha_i)},\]
		and
		\[\left|\prod^d_{i=1}\lambda_i\right|=\left|\dfrac{\mathrm{Res}(g,f)}{\mathrm{Disc}(f)}\right|,\]
		as was required to show.
	\end{proof}
\end{proposition}

\begin{definition}
	Let
	\begin{align*}
	\psi_{d}:\mathbb{R}&\longrightarrow \mathbb{R}^{d}\\
	r&\mapsto (r^{d-1}\ldots,r,1)^t
	\end{align*}
	define a curve in dimension $d$. For $n$ real numbers $r_1>r_2>\ldots>r_n$ we define the \textit{cyclic polytope} $\mathcal{C}(r_1,\cdots,r_n)$ to be the convex hull of those points on the curve $\psi_d$ corresponding to the $r_i$, i.e. \[
	\mathcal{C}(r_1,\cdots,r_n)=\mathrm{conv}\{\psi_{d}(r_i)|i=1,\ldots,n\}.
	\] 
\end{definition}

Let $f=\prod^d_{i=1}(x-\alpha_i)$, we write $\mathcal{C}(f)=\mathcal{C}(\alpha_1,\cdots,\alpha_d)\subset \mathbb{R}^{d}$ for a cyclic polytope of the roots of $f$ in dimension $d$. The associated matrix is $A_{\mathcal{C}(f)}$, where $(A_{\mathcal{C}(f)})_{ij}=\alpha_j^{d-i}$. It is known that $\det(A_{\mathcal{C}(f)})=\prod_{1\leq i<j\leq d}(\alpha_i-\alpha_j)$ (\cite[p. 11,]{25}). Furthermore, from the proof of Proposition \ref{prop:1-2}, we deduce that
\begin{align}\label{eq:3}
A_{\mathcal{K}(f)}^{-1}=DA_{\mathcal{C}(f)}^t,
\end{align}
where $D=\mathrm{diag}({f'(\alpha_1)^{-1}},\ldots,{f'(\alpha_d)^{-1}})$. For an irreducible polynomial $f$ we have
\[A_{\mathcal{C}(f)}A_{\mathcal{C}(f)}^t=(\mathrm{tr}_{K/\mathbb{Q}}(\alpha^{2d-i-j}))\in \mathrm{Mat}(d,\mathbb{Z}).\] 
In particular, for any $\beta\in R^{\vee}$, let $B=\mathrm{diag}(\sigma_1(\beta),\ldots,\sigma_d(\beta))$. Then  
\begin{equation}\label{eq:4}
A_{\mathcal{C}(f)}BA_{\mathcal{C}(f)}^t=(\mathrm{tr}_{K/\mathbb{Q}}(\beta\alpha^{2d-i-j}))\in\mathrm{Mat}(d,\mathbb{Z}).
\end{equation}

\begin{proposition}\label{prop:1-3}
	Let $f\in\mathbb{Z}[x]$ be an irreducible polynomial, and let $\alpha\in\mathbb{R}$ be a root of $f$. Then there exists a one-to-one correspondence between the interlacing polynomials of $f$ and elements in $\{\gamma\in\mathbb{Z}[\alpha]^{\vee}\mid\gamma\gg 0, \; \mathrm{tr}_{K/\mathbb{Q}}(\gamma)=1\}$.
	\begin{proof}
		Let us consider an injective homomorphism 
		\begin{align*}
		j:\mathbb{Z}[\alpha]^{\vee}&\longrightarrow\mathbb{R}^d\\
		\gamma&\mapsto (\sigma_1(\gamma),\ldots,\sigma_d(\gamma))^t.
		\end{align*}
		Given that $\mathbb{Z}[\alpha]^{\vee}=\dfrac{1}{f'(\alpha)}\mathbb{Z}[\alpha]$ (\cite[Prop. 2.2]{25}), we have (by Lemma \ref{lem:1-1}) a one-to-one map defined by $A_{\mathcal{K}(f)}^{-1}$:
		\begin{align}
		J:\mathbb{Z}^d&\longrightarrow j(\mathbb{Z}[\alpha]^{\vee}).\label{eq:5}
		\end{align} 
		By Proposition \ref{prop:1-1}, a polynomial $g\in\mathbb{Z}[x]$ interlaces $f$ if and only if there exists $\lambda_1,\ldots,\lambda_d\in \mathbb{R}_+$ such that $\sum\lambda_i=1$, and $g=\sum^d_{i=1}\lambda_i\prod^d_{j\not=i}(x-\alpha_j)$. The map (\ref{eq:1}) associates such a polynomial $g$ with the vector $b\in \mathbb{Z}^d$ such that $A_{\mathcal{K}(f)}\lambda=b$. Therefore $J(b)\in j(\mathbb{Z}[\alpha]^{\vee})$ corresponds to a totally positive element of trace one in $\mathbb{Z}[\alpha]^{\vee}$, as was required to show.
	\end{proof}
\end{proposition}

Let us notice that if we remove the trace requirement, then we would need to allow for non-monic polynomials that interlace $f$. Thus counting elements of trace $t$ is equivalent to counting integral points in $t\mathcal{K}(f)$. 

Before we move in to counting the interlacing polynomials, let us list some of the properties of sets $\mathcal{C}(f)$ and $\mathcal{K}(f)$.

\begin{corollary}\label{c:1-1}
	Let $f\in\mathbb{Z}[x]$ be a separable polynomial of degree $d$ such that all its roots are real. Then $|\Lambda_{\mathcal{C}(f)}|=|\mathcal{C}(f)\cap \mathbb{Z}^d|$.
	\begin{proof}
		We have $r\in \mathbb{Z}^d\cap\mathcal{C}(f)$ if and only if there exists $\lambda \in \Lambda_{\mathcal{C}(f)}$ such that $A_{\mathcal{C}(f)}\lambda=r$. Therefore it suffices to show that there cannot exist $\mu \in  \Lambda_{\mathcal{C}(f)}$ such that $\mu\not=\lambda$ and $A_{\mathcal{C}(f)}\mu=r$. Given that $f$ is a separable polynomial implies that $A_{\mathcal{C}(f)}$ is an invertible matrix, the corollary follows.
	\end{proof}
\end{corollary}

\begin{proposition}\label{prop:1-4}
	For a separable $f \in \mathbb{Z}[x]$ we have $\Lambda_{\mathcal{K}(f)}=\Lambda_{\mathcal{C}(f)}$.
	\begin{proof}
		Let $\lambda \in \Lambda_{\mathcal{K}(f)}$. Using identity (\ref{eq:3}) we have from $A_{\mathcal{K}(f)}\lambda=b$ that
		\begin{align*}
		\lambda&=DA_{\mathcal{C}(f)}^tb\\
		A_{\mathcal{C}(f)}\lambda&=A_{\mathcal{C}(f)}DA_{\mathcal{C}(f)}^tb,
		\end{align*}
		where both $b$ and $A_{\mathcal{C}(f)}DA_{\mathcal{C}(f)}^t$ are integral (see (\ref{eq:4})), thus $\lambda \in \Lambda_{\mathcal{C}(f)}$. The reverse inclusion is analogous.  
	\end{proof}
\end{proposition}

\begin{corollary}\label{c:1-2}
	Let $f\in \mathbb{Z}[x]$ be a separable monic polynomial such that all its roots are real. Then $|\mathcal{C}(f)\cap \mathbb{Z}^d|=|\mathcal{K}(f)\cap \mathbb{Z}^d|$.
\end{corollary}  

\begin{example}\label{ex:1-1}
	Let $f=x^2-D$, where $D \in \mathbb{N}$. Then 
	\begin{align*}
	\mathcal{K}(f)&=\{(1,\lambda\sqrt{D}-(1-\lambda)\sqrt{D})^t\mid \lambda \in [0,1]\}\\
	&=\{(1,(2\lambda-1)\sqrt{D})^t\mid \lambda \in [0,1]\}\\
	&=\{(1,\lambda'\sqrt{D})^t\mid \lambda' \in [-1,1]\}.
	\end{align*}	
	Therefore all the integer points in $\mathcal{K}(f)$ correspond to integers in the interval $[-\sqrt{D},\sqrt{D}]$, and from equation (\ref{eq:2}), the elements of $\Lambda_{\mathcal{C}(f)}$ correspond to $\left(\dfrac{-\sqrt{D}+r}{-2\sqrt{D}},\dfrac{\sqrt{D}+r}{2\sqrt{D}}\right)^t$ for $r\in [-\sqrt{D},\sqrt{D}]\cap\mathbb{Z}$. In particular, 
	\begin{equation}\label{eq:6}|\mathcal{K}(f)\cap \mathbb{Z}^2|=2\lfloor \sqrt{D}\rfloor+1\approx \sqrt{\mathrm{Disc}(f)}.
	\end{equation}
	A similar result holds for quadratic polynomials of the form $x^2-x-D$.
\end{example}

\begin{corollary}\label{c:1-3}
	The number of interlacing polynomials of an irreducible quadratic polynomial goes to infinity with the discriminant of the polynomial. 
\end{corollary}

We would like to replicate the above result for polynomials of higher degrees. More generally, we wish to determine the conditions under which a given polynomial is interlaced, and count the interlacing polynomials. We cannot apply Minkowski's First Theorem (\cite{25}), as the associated convex set to a polynomial is not symmetric and generally, it does not contain the origin. We shall use the following variant of the Flatness Theorem:

\begin{theorem}\textup{\cite[Remark 2.7]{1}}\label{thm:1-6}
	Let $S$ be a simplex in $\mathbb{R}^d$. Then 
	\[w(S,\mathbb{Z}^d)\leq cd(1+\log(d)+|S\cap\mathbb{Z}^d|^{1/d}),\]
	where $c$ is a universal constant.
\end{theorem}

Observe that $\mathcal{C}(f)$ (and $\mathcal{K}(f)$) can be projected onto a simplex $S$ in $\mathbb{R}^{d-1}$. In particular, if $\mathcal{C}(f)=\mathrm{conv}(b_1,\ldots,b_d)$, then $S=\mathrm{conv}(b_1',\ldots,b'_d)$, where $(b'_i)_j=(b_i)_{j}$ for $1\le j\le d-1$. Therefore, without loss of generality, we can consider $\mathcal{C}(f)$ (and $\mathcal{K}(f)$) as a simplex. 

\begin{proposition}\label{prop:1-5}
	Let $f\in \mathbb{Z}[x]$ be a monic irreducible polynomial such that totally real algebraic integer $\alpha$ is a root of $f$. Then 
	\[
	w(\mathcal{C}(f),\mathbb{Z}^d)=\min\{\mathrm{Span}(\gamma)\mid \gamma\in\mathbb{Z}[\alpha]\setminus\mathbb{Z}\}.
	\]
	
	\begin{proof}
		Let $(\mathbb{Z}^d,\langle\,,\rangle)$ denote the integer lattice of rank $d$, where $\langle\,,\rangle$ is the usual inner product, i.e. $\langle v,w\rangle=\sum^d_{i=1}v_iw_i$ for $v,w\in \mathbb{R}^d$. Let $A$ be the associated matrix of $\mathcal{C}(f)$. Thus every element ${\gamma}\in \mathcal{C}(f)$ may be written as ${\gamma}=A{\lambda}$, where ${\lambda}\in \Lambda=\{{\theta}\in \mathbb{R}^d_+\mid \sum^d_{i=1}\theta_i=1\}$. Let $v\in \mathbb{Z}^d$, then 
		\begin{align*}
		\max_{w \in \mathcal{C}(f)}\langle v,w\rangle=	\max_{\lambda \in \Lambda}v^tA{\lambda}.
		\end{align*}
		As $v^tA=(\sum^{d}_{i=1}v_i\alpha^{d-i}_1,\ldots,\sum^{d}_{i=1}v_i\alpha^{d-i}_d),$ it follows that
		\begin{align*}
		\max_{w \in \mathcal{C}(f)}\langle v,w\rangle=	\max_j \sum^{d}_{i=1}v_i\alpha^{d-i}_j.
		\end{align*} 
		Similarly we have that 
		\begin{align*}
		\min_{w \in \mathcal{C}(f)}\langle v,w\rangle =\min_j \sum^{d}_{i=1}v_i\alpha^{d-i}_j.
		\end{align*} 
		For $v\in\mathbb{Z}^d$ let $g\in \mathbb{Z}[x]$ such that $g=\sum^{d}_{i=1}v_ix^{d-i}$. Therefore 
		\[
		\max_{w \in \mathcal{C}(f)}\langle v,w\rangle-\min_{w \in \mathcal{C}(f)}\langle v,w\rangle=\mathrm{Span}(\gamma),
		\]
		where $\gamma=g(\alpha)\in \mathbb{Z}[\alpha]$. By the definition of the width the proposition follows.
	\end{proof}
\end{proposition}

For a totally real algebraic integer $\beta \in \mathbb{R}$ there can exist infinitely many (up to equivalence) totally real algebraic integers $\alpha$ of a bounded degree such that $\beta \in \mathbb{Z}[\alpha]$. For example, consider the family of polynomials 
\[
\{x^4-2ax^2+a^2-2\in \mathbb{Z}[x]\mid a\ge2\}.
\]
The roots of these polynomials are $\pm\sqrt{a\pm\sqrt{2}}$. Thus for each of the corresponding simplices, the width is bounded above by $2\sqrt{2}$. More generally, it suffices to consider number fields of bounded degree that contain $\mathbb{Q}[\beta]$ as a subfield.

\theoremstyle{plain}
\newtheorem*{thm:1}{Theorem \ref{Thm:1}}
\begin{thm:1}
	Let $d,r \in \mathbb{N}$. Up to $\mathbb{Z}$-equivalence, there are only finitely many irreducible monic polynomials $f\in \mathbb{Z}[x]$ of degree $d$ such that $\mathbb{Q}[x]/(f)$ is a primitive field and $f$ is interlaced by at most $r$ monic integer polynomials.
	\begin{proof}
		
		First, we claim that, up to $\mathbb{Z}$-equivalence, there are only finitely many $g\in\mathbb{Z}[x]$ of degree $d$ of bounded span. Without loss of generality, we assume that all roots of $g$ are positive, and at least one of the roots is less than one. Given that $|\{\sum^d_{i=0}a_ix^i\mid |a_j|<M\}|$ is a finite set, our claim follows. 
		
		Let $f\in \mathbb{Z}[x]$ such that $K=\mathbb{Q}[x]/(f)$ is a primitive field, and $\alpha\in\mathbb{R}$ such that $f(\alpha)=0$. If $\beta \in\{\gamma\in\mathbb{Z}[\alpha]\mid \mathrm{span}(\gamma)=\min_{ \delta\in \mathbb{Z}[\alpha]}\mathrm{span}(\delta)\}$, then there are only finitely many algebraic integers $\gamma\in K$ (up to $\mathbb{Z}$-equivalence) such that $\beta \in \mathbb{Z}[\gamma]$ (Corollary 6.2.2. in \cite{11}). From Proposition \ref{prop:1-5}, it follows that up to $\mathbb{Z}$-equivalence, there are only finitely many such $f$ with bounded width, and in the light of Theorem \ref{thm:1-6}, follows the theorem.  
	\end{proof}
\end{thm:1}

\section{Non-interlacing polynomials}
Here we demonstrate that there exist irreducible polynomials that are non-interlacing, but, not for degrees $2,3,4,5$ and $7$. The following proposition was proved by Dobrowolski for minimal polynomials of integer symmetric matrices. We shall reprove it, using the results from the previous section:

\begin{proposition}\label{prop:1-6}\textup{\cite[Lemma 1]{8}}
	Let $f\in \mathbb{Z}[x]$ be a monic and irreducible polynomial of degree $d$. If $f$ is interlaced then $|\mathrm{Disc}(f)|\ge d^d$.
	\begin{proof}
		Let $f=\prod^d_{i=1}(x-\alpha_i)\in\mathbb{Z}[x]$ be interlaced by $g$. Then there exists $\lambda \in \Lambda_{\mathcal{C}(f)}$ such that $g=\sum^d_{i=1}\lambda_i\prod_{j\not= i}(x-\alpha_j)\in\mathbb{Z}[x]$. By Proposition \ref{prop:1-2} it follows that 
		\begin{align*}
		\left|\dfrac{\mathrm{Res}(f,g)}{\mathrm{Disc}(f)}\right|&=\left|\prod^d_{i=1}\lambda_i\right|\\
		&\le\left|\dfrac{1}{d} \sum^d_{i=1}\lambda_i\right|^d,
		\end{align*}
		and so
		\[\left|d^d\mathrm{Res}(f,g)\right|\le\left|\mathrm{Disc}(f)\right|.\]
		
		As $f,g\in\mathbb{Z}[x]$ are monic polynomials, we conclude that $\mathrm{Res}(f,g)\in\mathbb{Z}\setminus\{0\}$, and the proposition follows.
	\end{proof}
\end{proposition}

There exist polynomials that satisfy the hypothesis of the theorem above, but have the discriminant smaller than $d^d$ (\cite{29}). The smallest such example known has degree 2880. We shall show that the smallest degree for which there exists an irreducible non-interlacing polynomial is 6. We begin by characterising an infinite family of non-interlacing polynomials. For this we need to use the following result: 

\begin{lemma}\label{lem:1-2}
	Let $M=\mathbb{Q}[\zeta_d^{}]$ such that $d$ is square-free, $d$ is not a prime number or twice a prime number, $d> 20$ and $d\not=30$. Let $S$ be the ring of integers in $M$. Let $K$ be the maximal totally real subfield of $M$, i.e. $K=\mathbb{Q}[\zeta_d^{}+\zeta^{-1}_d]$, with its ring of integers $R$. Then the minimum of any ideal $R$-lattice is at least 2.
	\begin{proof}
		Let $I$ be an ideal in $R$, and let $\gamma \in K$ such that $(I,\mathrm{tr}_{K/\mathbb{Q}}(\gamma xy))$ is an ideal lattice. Let $m=\min(I,\mathrm{tr}_{K/\mathbb{Q}}(\gamma xy))$ and $z \in \mathcal{M}(I,\mathrm{tr}_{K/\mathbb{Q}}(\gamma xy))$. We extend this lattice to $S$, let 
		\[I^e=\{\sum a_ib_i|a_i \in I,\;b_i \in S\},\] 
		and $\mathrm{tr}_{M/K}(\mathrm{tr}_{K/\mathbb{Q}}(\gamma xy))= \mathrm{tr}_{M/\mathbb{Q}}(\gamma x\overline{y})$. Clearly $(I^e,\mathrm{tr}_{M/\mathbb{Q}}(\gamma x\overline{y}))$ is an ideal $S$-lattice. Given that $z \in I^e$, we have that the minimum of $\mathrm{tr}_{M/\mathbb{Q}}(\gamma x\overline{y})$ is smaller than $\mathrm{tr}_{M/\mathbb{Q}}(\gamma z^2)$. Thus
		\begin{align*}
		\mathrm{tr}_{M/\mathbb{Q}}(\gamma z^2)&=\mathrm{tr}_{M/K}(m)\\
		&=2m\\
		&\ge 4,
		\end{align*}
		the last inequality follows from \cite[Lemma 1.4 \& Cor. 2.2]{2}, thus $m\ge 2$.
	\end{proof}
\end{lemma}
From the lemma above and Proposition \ref{prop:1-3} follows:

\theoremstyle{plain}
\newtheorem*{thm:2}{Theorem \ref{Thm:2}}
\begin{thm:2}
	Let $d\in\mathbb{N}$ such that $d$ is square-free, $d$ is not a prime number or twice a prime number, $d> 20$ and $d\not=30$. Then the minimal polynomial of $\zeta_{d}^{}+\zeta_{d}^{-1}$ is non-interlacing. 
\end{thm:2} 

An example of such polynomial of smallest degree is the minimal polynomial of $\zeta_{21}^{}+\zeta_{21}^{-1},$ 
\begin{equation}\label{eq:7}
x^6 - x^5 - 6x^4 + 6x^3 + 8x^2 - 8x + 1.
\end{equation}
It is known that for primes $p$ larger than $5$, the minimal polynomial of $\zeta_p^{}+\zeta_p^{-1}$ is always interlaced (\cite{13}). A classification of interlacing minimal polynomials of $\zeta_d^{}+\zeta_d^{-1}$ in relation to integer symmetric matrices can be found in \cite{24}. 

Let us define $\sigma_r(I)=\sum_{J\mid I}\mathrm{N}(J)^r,$ where $I$ is an ideal of $R$. Furthermore, let 
\[s^{K}_{l}(m)=\sum_{\substack{\gamma \in R^{\vee}_+\\ \mathrm{tr}_{K/\mathbb{Q}}(\gamma)=l}}\sigma_{m-1}((\gamma)(R^{\vee})^{-1}),\]
where $(R^{\vee})^{-1}$ is the \textit{different} ideal.

\begin{theorem}\textup{\cite{27,31}}
	Let $K$ be a totally real algebraic number field of degree $d>1$ and let $\zeta_K$ be its Dedekind zeta function. Let $h=2dm$, where $m \in \mathbb{N}$. Then
	\[
	\zeta_{K}(1-2m)=2^d\sum_{l=1}^r b_l(h)s^{K}_l(2m)\in \mathbb{Q}\setminus\{0\}.
	\] 
	The numbers $r\ge 1, b_1(h),\ldots,b_r(h)$ are rational and they only depend on $h$, where
	\begin{align*}
	r=\begin{cases} \left\lfloor \dfrac{h}{12}\right\rfloor & \mbox{ if } h\equiv 2 \pmod{12}\\\\
	\left\lfloor \dfrac{h}{12}\right\rfloor+1 & \mbox{ if } h\not\equiv 2 \pmod{12}.
	\end{cases}\\[-\normalbaselineskip]\tag*{}
	\end{align*}
\end{theorem}

\theoremstyle{plain}
\newtheorem*{thm:3}{Theorem \ref{Thm:3}}
\begin{thm:3}
	There does not exist a non-interlacing irreducible integer polynomial of degree $2,3,4,5$ or $7$.
	\begin{proof}
		Let $f \in \mathbb{Z}[x]$ be an irreducible polynomial of degree $2,3,4,5$ or $7$, such that $K\cong\mathbb{Q}[x]/(f)$ with the ring of integers $R$,  and $\alpha\in\mathbb{R}$ be a root of $f$. Given that $\mathbb{Z}[\alpha]\subseteq R$ implies that  $R^{\vee}\subseteq\mathbb{Z}[\alpha]^{\vee}$. From the theorem above it follows that for $d=2,3,4,5$ or $7$, we have 
		\[\zeta_{K}(-1)=2^d b_1(h)\sum_{\substack{\gamma \in R^{\vee}_+\\ \mathrm{tr}_{K/\mathbb{Q}}(\gamma)=1}}\sigma_{1}((\gamma)(R^{\vee})^{-1}).\]
		As $\zeta_{K}(-1)\neq 0$, therefore there exists at least one $\gamma\in R^{\vee}_+\subseteq\mathbb{Z}[x]^{\vee}$ such that $\mathrm{tr}_{K/\mathbb{Q}}(\gamma)=1$. By Proposition \ref{prop:1-3} follows the proof.
	\end{proof}
\end{thm:3}

The requirement for irreducibility is necessary. For example, $(x-1)(x-2)$ is non-interlacing. 

\begin{theorem}\textup{\cite{12}}
	Let $K$ be a totally real algebraic number field of degree $d$ and let $R$ be its ring of integers. For $1\le d \le 7$, if $K$ satisfies Condition (A), then there always exists an element $\alpha\in R^{\vee}_{+}$ such that $\mathrm{tr}_{K/\mathbb{Q}}(\alpha)=1$.
	\begin{proof}
		Let $f$ be the minimal polynomial of $\alpha$ such that $R=\mathbb{Z}[\alpha]$. From Theorem \ref{Thm:3} it follows that for $d=2,3,4,5,$ and $7$, polynomial $f$ is interlacing. This leaves us to prove the theorem for sextic number fields. Given that there exist units of every signature in $K$, implies that the codifferent is narrowly equivalent to a square of an ideal in $K$ (\cite[Theorem 176]{15}), i.e. there exists a fractional ideal $I$ in $K$, such that $\gamma I^2=R^{\vee}$, where $ \gamma\in K$ and $\gamma\gg 0$. Therefore,  $(I,\mathrm{tr}_{K/\mathbb{Q}}(\gamma xy))$ is a positive definite unimodular $\mathbb{Z}$-lattice.  From the classification of unimodular quadratic forms (\cite[106:13]{26}) it follows that there exists an element of trace one in $R^{\vee}_+$. 
	\end{proof}
\end{theorem}

The number field $\mathbb{Q}(\zeta_{21}^{}+\zeta_{21}^{-1})$ is monogenic (see \cite[Prop. 2.16]{30}), and given that the polynomial (\ref{eq:7}) is non-interlacing, this implies that there cannot exist units of every signature in $\mathbb{Z}[\zeta_{21}^{}+\zeta_{21}^{-1}]$.

\section{Universal quadratic forms}

For a quadratic form $Q$ of rank $r$ over $R$, in what follows, we shall write $L$ for $R^r$. We begin by proving the following useful proposition:

\begin{proposition}\label{prop:1-7}
	Let $K$ be a totally real algebraic number field satisfying Condition (A), and let $R$ be its ring of integers. Let $Q$ be a universal quadratic form over $R$. Then there exists $\delta \in K_+$ and a positive definite $\mathbb{Z}$-lattice $(L,\mathrm{tr}_{K/\mathbb{Q}}(B_Q^{\delta}))$ such that 
	\begin{equation}\label{eq:8}
	|\mathcal{M}(L,\mathrm{tr}_{K/\mathbb{Q}}(B^{\delta}_Q))|\ge 2|\mathcal{M}(K)|.
	\end{equation}
	
	\begin{proof}
		Given that $K$ satisfies Condition (A), implies that there exists $\delta\in K_+$ such that $R^{\vee}=\delta R$. Therefore, $Q^{\delta}$ represents all the elements in $R^{\vee}_+$. Let us consider a $\mathbb{Z}$-lattice $(L,\mathrm{tr}_{K/\mathbb{Q}}(B_Q^{\delta}))$ (replacing $\delta$ by $2\delta$ if necessary). We note that $\min(L,\mathrm{tr}_{K/\mathbb{Q}}(B_Q^{\delta}))$ is equal to $\min_{\gamma \in R^{\vee}_+}\mathrm{tr}_{K/\mathbb{Q}}(\gamma )$ (or twice as much). Thus
		\begin{align*}
		|\mathcal{M}(L,\mathrm{tr}_{K/\mathbb{Q}}(B_Q^{\delta}))|&\ge 2|\mathcal{M}(K)|,
		\end{align*}
		where a factor of two on the right follows from the fact that $Q(\pm v)=\alpha$ for all $\alpha\in R_{+}$.
	\end{proof}
\end{proposition}

\begin{example}
	Let $Q$ be a quadratic form 
	\[
	x^2+2y^2+2yz+(2+\epsilon)z^2
	\]
	over $\mathbb{Z}[\frac{1+\sqrt{5}}{2}]$, where $\epsilon=\frac{1+\sqrt{5}}{2}$. From Theorem 1.1 in \cite{7} we know that $Q$ is a universal quadratic form. Given that $\epsilon$ is a unit of negative norm implies that $\mathbb{Q}(\sqrt{5})$ satisfies Condition (A), and thus fulfils the hypothesis of the above proposition. The minimal polynomial of $\frac{1+\sqrt{5}}{2}$ is $f=x^2-x-1$. It is interlaced only by two polynomials, $x$ and $x-1$. By Proposition \ref{prop:1-3}, it follows that $|\mathcal{M}(L,\mathrm{tr}_{\mathbb{Q}(\sqrt{5})/\mathbb{Q}}(B_Q^{\delta}))|\ge4$ (note that this bound will hold for all the lattices over $\mathbb{Z}[\frac{1+\sqrt{5}}{2}]$, invariant of $Q$). On the other hand, let 
	\[
	\delta=\dfrac{\epsilon}{f'(\epsilon)}=\dfrac{\epsilon+2}{5},
	\]
	and by computing the minimal vectors of the $\mathbb{Z}$-lattice we conclude that (\ref{eq:8}) is an equality in this case. However, we can have a strict inequality also. Let us consider an another universal quadratic form $Q'$ over $\mathbb{Z}[\frac{1+\sqrt{5}}{2}]$:
	\[x^2+y^2+z^2.\]
	In this case $|\mathcal{M}(L,\mathrm{tr}_{\mathbb{Q}(\sqrt{5})/\mathbb{Q}}(B_{Q'}^{\delta}))|=12$.
\end{example}

\begin{corollary}\label{c:1-4}
	Let $K$ be a totally real algebraic number field of degree $d$ satisfying Condition (A), and let $R$ be its ring of integers. Let $r$ be the rank of a universal quadratic form over $R$. Then
	\[
	r\ge \dfrac{\log(2|\mathcal{M}(K)|)}{0.401d(1+o(1))\log(2)}.
	\]
	\begin{proof}
		As in Proposition \ref{prop:1-7}, a given quadratic form of rank $r$ over the ring of integers of $K$ can be lifted to a rational integer lattice of rank $dr$. From the bound on the kissing number $\tau_{dr}\le 2^{0.401dr(1+o(1))}$ (\cite{17}) and inequality (\ref{eq:8}) follows the corollary.
	\end{proof}
\end{corollary}

This bound can be improved to
\[
r\ge \dfrac{2|\mathcal{M}(K)|}{\tau_d}
\]
for universal diagonal forms. As a special case, we have the following theorem:

\begin{theorem}\label{thm:1-7}
	Let $K$ be a totally real algebraic number field of degree $d$ satisfying Condition (A), and let $R$ be its ring of integers. If there exists an element $\alpha\in R^{\vee}_{+}$ such that $\mathrm{tr}_{K/\mathbb{Q}}(\alpha)=1$, then the rank of a universal quadratic form $Q$ over $R$ is larger than
	\[\begin{cases}
	\dfrac{|\mathcal{M}(K)|}{d} & \mbox{if $Q$ is a classical quadratic form}\\
	\dfrac{1+\sqrt{1+4|\mathcal{M}(K)|}}{2d} & \mbox{if $d>5$ or $|\mathcal{M}(K)|>240$}\\
	\dfrac{|\mathcal{M}(K)|}{15d}& \mbox{otherwise}.
	\end{cases}\]
	\begin{proof}
		Observe that for a classical quadratic form if $\min_{v \in L}(L,B_Q)=1$, then $L\cong <1>\perp L'$. Therefore, if $L$ is a $\mathbb{Z}$-lattice, then the rank of $L$ is larger or equal to $\mathcal{M}(L)/2$. Consider a $\mathbb{Z}$-lattice $(L,\mathrm{tr}_{K/\mathbb{Q}}(B_Q^{\delta}))$ as in Proposition \ref{prop:1-7}, where $Q$ is a classical quadratic form of rank $r$, and thus the rank of $L$ over $\mathbb{Z}$ is $dr$. By Proposition \ref{prop:1-7} we have 
		\begin{align*}
		\dfrac{|\mathcal{M}(L)|}{2}&\ge |\mathcal{M}(K)|\\
		dr&\ge |\mathcal{M}(K)|\\
		r&\ge \dfrac{1}{d}|\mathcal{M}(K)|.
		\end{align*}
		On the other hand, if $Q$ is an integral quadratic form, then $2Q$ is a classical quadratic form. The minimum of $\mathbb{Z}$-lattice $(L,\mathrm{tr}_{K/\mathbb{Q}}(B_Q^{2\delta}))$ is 2. Let $M$ be a lattice generated by the minimum vectors of $L$. It follows that $M$ is a direct sum of root lattices (see Theorem 4.10.6 in \cite{23}) and $|\mathcal{M}(M)|=|\mathcal{M}(L)|$. Considering root lattices $A_n$, $D_n$ and $E_8$, we observe that $|\mathcal{M}(A_n)|\le |\mathcal{M}(D_n)|$ for $n\ge 5$ (see Table 4.10.13 in \cite{23}). Given that $|\mathcal{M}(D_{n})|=2n(n-1)$ and $|\mathcal{M}(E_8)|=240$, if $dr> 16$, we can bound the number of minimal vectors in $L$ by $2dr(dr-1)$. For the case when $dr\ge 16$ or $|\mathcal{M}(K)|\le 240$, we can bound the number of minimal vectors by $30dr$. Given that there cannot exist a positive definite binary universal quadratic form, we assume that $r\ge 3$. The rest of the proof follows as for the classical quadratic form. 
	\end{proof}
\end{theorem}

\begin{corollary}\label{c:1-5}
	The minimal rank of a universal quadratic form over quadratic fields that contains a unit of negative norm grows at the order of the fourth root of the discriminant.
	\begin{proof}
		This follows from the theorem above and equation (\ref{eq:6}).
	\end{proof}
\end{corollary}

The above result appeared previously as Proposition 5 in \cite{6}, however, here it is independent of the Riemann hypothesis. In the light of Theorem \ref{Thm:1} we have:

\theoremstyle{plain}
\newtheorem*{thm:4}{Theorem \ref{Thm:4}}
\begin{thm:4}
	Let $d,r \in \mathbb{N}$. There are only finitely many totally real primitive number fields of degree $d$ that satisfy Condition (A) and have an integer universal quadratic form of rank $r$ defined over them.  
	\begin{proof}
		Let $f$ be the minimal polynomial of $\alpha$ such that $R=\mathbb{Z}[\alpha]$. From the proof of Theorem \ref{Thm:1}, we know that there are only finitely many fields for which there are at most $r$ interlacing polynomials. Therefore, the number of interlacing polynomials grows, and from Proposition \ref{prop:1-3}, we deduce that the right side of inequality (\ref{eq:8}) grows also. Given that the maximal number of the minimal vectors is bounded by the kissing number, follows that the rank of a universal quadratic form over such number fields should grow also, as was required to show. 
	\end{proof}
\end{thm:4} 

Quadratic number fields are always monogenic, thus the main difficulty in satisfying Condition (A) is to show that there exists a unit of negative norm. For number fields of higher degrees, along with proving that the field has units of every signature, we will have to prove that the field is monogenic. In general, it is unclear how to meet those two conditions. But, for a given field of low degree, it can be readily checked. Let us consider an example:

\begin{example}
	Let $f=x^5 - 2x^4 - 9x^3 + 4x^2 + 15x + 3$. It is an irreducible polynomial with only real roots. Let $K$ be the number field generated by the root of $f$, let $R$ be its ring of integers. We note that $\mathrm{Disc}(f)=3\times61\times241\times547$, therefore $K$ is monogenic, in particular, $R=\mathbb{Z}[\alpha],$ where $\alpha$ is a root of $f$. Furthermore, the narrow class number of $R$ is $1$, thus we have units of every signature, and therefore, $K$ satisfies Condition (A). We computed that $f$ has $218$ interlacing polynomials, thus if there exists a classical universal quadratic form over $K$, then it has to be at least of rank $44$.
\end{example} 

Dummit and Kisilevsky,  in \cite{9}, studied a parametric family of cubic polynomials
\begin{equation*}\label{pol:1-1}
f_k=x^3+x^2-(3k^2 + k + 2)x - (2k^3 + 2k^2 + 2k + 1),
\end{equation*}
related to the cubic subfield of $\mathbb{Q}(\zeta_d^{})$. They proved that for infinitely many $k$, a root of $f_k$ forms a power integer basis of the ring of integers.

\theoremstyle{plain}
\newtheorem*{thm:5}{Theorem \ref{Thm:5}}
\begin{thm:5}
	For any given $r\in \mathbb{Z}$, there exist infinitely many totally real quadratic and cubic number fields that do not have a universal quadratic form of rank $r$ defined over them.
	\begin{proof}
		This clearly holds for real quadratic number fields, as it a classical result about Pell equations that there are infinitely many real quadratic number fields that have a unit of negative norm.
		
		For cubic number fields, let $\alpha\in \mathbb{R}$ be a root of $f_k$, let $K\cong\mathbb{Q}[x]/(f_k)$ be a number field with the roots of integers $R$. As for infinitely many $k$, the ring $K$ is monogenic (\cite[Theorem 3]{9}), it suffices to show that a proportion of them has units of every signature. This will give us Condition (A), and the proof of this theorem will follow from Theorem \ref{Thm:4}.
		
		Both $\alpha+k$ and $\alpha+k+1$ are units in $R$, as $f_k(-k)=-1$ and $f_k(-k-1)=1$. We claim that $\alpha+k, \alpha+k+1$ and $-1$ generate units of every signature.
		This follows from the fact that the $\mathrm{N}_{K/\mathbb{Q}}(\alpha+k)=1$ and $\mathrm{N}_{K/\mathbb{Q}}(\alpha+k+1)=-1$, thus $\alpha+k$ is not totally positive. Letting $\epsilon_1=\alpha+k$ and $\epsilon_2=\alpha+k+1$, then $\{\pm1,\pm\epsilon_1,\pm\epsilon_2,\pm\epsilon_1\epsilon_2\}$ contains a unit of every possible signature, as was claimed.
	\end{proof}
\end{thm:5}

\end{document}